# STRUCTURAL DYNAMICS OF VARIOUS CAUSES OF MIGRATION IN JAIPUR

Dr. Jayant Singh, Assistant Professor, Department of Statistics,

University of Rajasthan, Jaipur, India

E-mail: jayantsingh47@rediffmail.com

Hansraj Yadav, Research Scholar, Department of Statistics,

University of Rajasthan, Jaipur, India

E-mail: hansraj_yadav@rediff.com

Dr. Florentin Smarandache, Chair of Math & Sciences Department,

University of New Mexico, Gallup, USA

E-mail: smarand@unm.edu

**Abstract**

Jaipur urban area has grown tremendously in last three decades. Composition of People migrating due to various reasons has display a meticulous trend. Dominance of people moving due to marriages is getting sturdy whereas Jaipur city is losing its luster in attracting persons for education and business. Short duration migration from Jaipur district to urban area has gone down to a very low level. Flow of migrants from Rural areas to Jaipur outpaced the migrants from urban areas and its composition from various in terms long and short distances migration has substantially changed over two consecutive decades. Movements of males and females were differ on many criterion as male moving faster than females for employment & education and females move faster than male for marriages and moving along family was found evident in short, medium and long distances migration. Gender gap in people migration from different reasons was observed and a gender specific trend was seen favour. Short duration migration and migration due to education & employment is not as prominence as it was two decade back

**INTRODUCTION:**

Migration from one area to another in search of improved livelihoods is a key feature of human history. While some regions and sectors fall behind in their capacity to support populations, others move ahead and people migrate to access these emerging opportunities. There are various causes like political, cultural, social, personal and natural forces but aspire for betterment, higher earning, more employment opportunities receive special attention. There are four type of migration namely

i. Rural-Rural

ii. Rural-Urban

iii. Urban- Urban

iv. Urban-Rural

Though all of these have different implication over the various demographic and socio-economic characteristics of the society but rural-urban & urban-urban migration is a cause of concern in reference to migration process to Jaipur urban agglomeration. The dynamics of migration for three census (1981, 1991, 2001) has been analyzed from different angles at destination i.e. Jaipur Urban Agglomeration. The peoples of two places have different socio-economic character like education attainment, availability of land to the rural labor and agriculture production capacity, industrialization etc and the difference of these factors at two places gear the migration process.

Distance plays a prominent role in migration of peoples, in general people from nearby area show a faster pace than the distant places due to psychic of being come back or feel like at home or the reason that some acquaintance in nearby area plays a big pull factor. However these assumptions do govern by other consideration of pull and push factor and the prevalent socio-economic aspects of the origin and destination places.

Jaipur being the capital of the state and proximity to the national state has been a great potential to draw peoples. It has not been attracting peoples from the nearby areas but it has influence on the persons of entire state and other states of the country. Majority of immigrants to Jaipur belongs to different parts of the states followed by its adjoining states. However it has been able to attract people from all over the country and overseas as well though their contribution in totality is not as significant. Seeing at this scenario it is worthwhile to limit the migrants from the following area to comprehend the migrant process of Jaipur. In-migrants to Jaipur urban area from (a) various parts of Jaipur district (b) other districts of the state (c) adjoining states of the state having fair share in migrants and (d) total migration which is overall migration from all the areas.

**COMPOSITION OF IN-MIGRANTS TO JAIPUR:**

In-migrants to Jaipur has grown by leaps and bounds in the last three decades. The decadal growth of in-migrants to Jaipur in last four decades synchronized with the growth of urban population of the Jaipur. Though the decadal rate of growth of migrants is lagging behind to the growth of the urban population as both has been 59.3% & 45.2% in decade 1991-2001, 49.5% & 35.8 % in decade 1991-81 respectively. Short Distance migration is considered, people from the other parts of Jaipur district who are coming to jaipur urban area, migration from other parts of the state is relatively longer distance migration and put in the moderate (medium) distance migration whereas the people from out side the state are in the category of long distance migration. The contribution of the short distance migration in total migration as per census of 2001, it was 17.1% against the 51% were medium distant migrants as they came from other districts of the state and long distance migration from some most contributing states namely Punjab, U. P., & Delhi have there share as 9.6%, 3.3% & 2.3% in total migrants to Jaipur in this same duration. These three states accounted for half of the long distance migration.

These different types of migration spell a meticulous trend over the years. As small distance migration shows a downward trend as its share in total migration which was 28.8% in yr 1981 came to 25.8% in according to census of 1991 and further slipped to 17.1% in census 2001. Medium distance migration exhibited a opposite path to the short distance migration as it advanced to 47% in yr 1991 against 45% in yr 1981 which further ascended to 51% in yr 2001. Contribution of long distance migration in total migration from all states also exhibited rolling down trend. This trend followed suite for the migration from the adjoining states.

**COMPOSITIONAL DYNAMICS OF REASON FOR MIGRATION TO JAIPUR URBAN AREA:**

Affect of various reasons of migration on peoples of diverse areas is different. Some reasons are more common than others moreover their affect on male and females is also different. Share of Rural and Urban in-migrants population will widely vary for various cause of migration. Distance of place of origin is also a crucial factor in migration process to any area. Dynamics of various reasons for migration will be analyzed from four perspectives.

1. Dominance of various reasons for migration;
2. Rural-Urban Paradigm and changes taking place;
3. Gender issues and disparity.

## DOMINANCE OF VARIOUS REASONS FOR MIGRATION IN MIGRATION PROCESS:

Person do migrate from a variety of reasons, prominent of them are migration due to 1. Employment 2. Education 3. Marriage 4. Moving with family. Marriage has been the foremost reason for migration as its share in total migrants to Jaipur was 32.1% in yr 2001. People migrating for the employment and/or business with 27.3% contribution in total migration seconded the marriage cause. It was distantly followed by category of persons moving with family with 17.6 % share in total migration. There was a remarkable difference in two dominating categories of people moving due to employment and marriage and it was that the people migrating to Jaipur due to employment is on declining side as it came down to 25% in yr 2001 from 27.3% in yr 1991 and 30.2% in yr 1981 contrary to a gradual increase in people migrating to Jaipur because of marriage as it raised to 32.1% in yr 2001 from 27.8% in yr 1991 and 25.2% in yr 1981.

Education as a cause of migration doesn't have significant contribution in total migration to Jaipur and it is getting meager over the years. As in yr 1981 its share in total migration was 6.1% and the figure came to 4.4% in yr 1991 and further dip to 2.7% in yr 2001. This movement is also followed by migrants for education from all the adjoining state, within state and from Jaipur district to Jaipur urban area. People moving with household also followed the decline suite though the rate of decline was steeper than the others as the share of people migrating under this category which was 30.2% in yr 1981 fall to 28.5% in yr 1991 and further it slip to 17.6 % in yr 2001. Composition of various reasons for migration over last three decades is depicted in coming Graph.

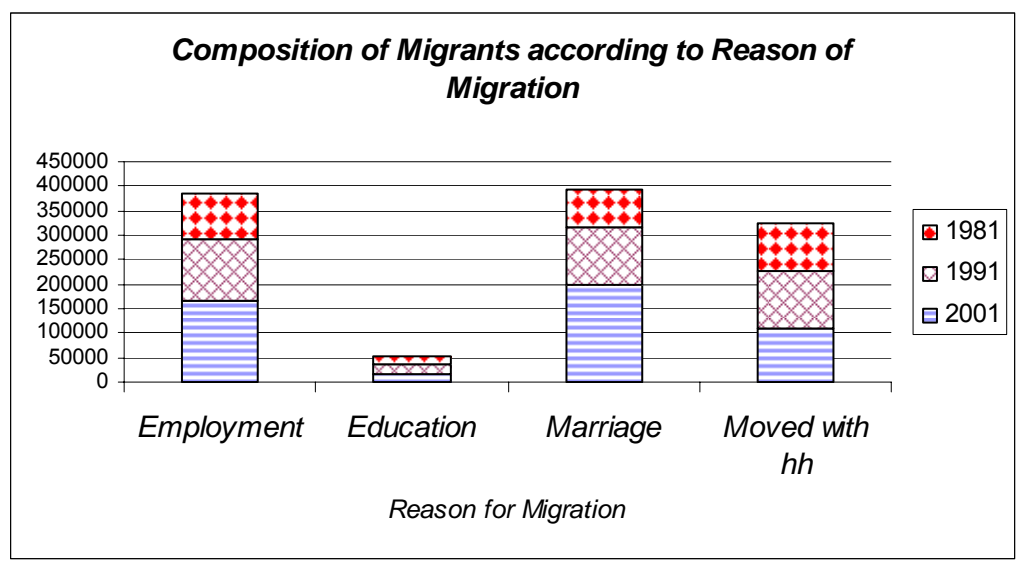

**RURAL-URBAN PARADIGM:**

Intensity of migration widely differs for persons migrating form Rural and Urban areas for various reasons for migration. Flow of migrants from Rural areas to Jaipur outpaced the migrants from urban areas. According to data of census in yr 1981, the share of migrants to Jaipur urban area from rural and urban areas was 53% & 47 % respectively and this gap remained intact in the coming decades. The trend in rural, urban and combined for last three decades is depicted in graph on next page.

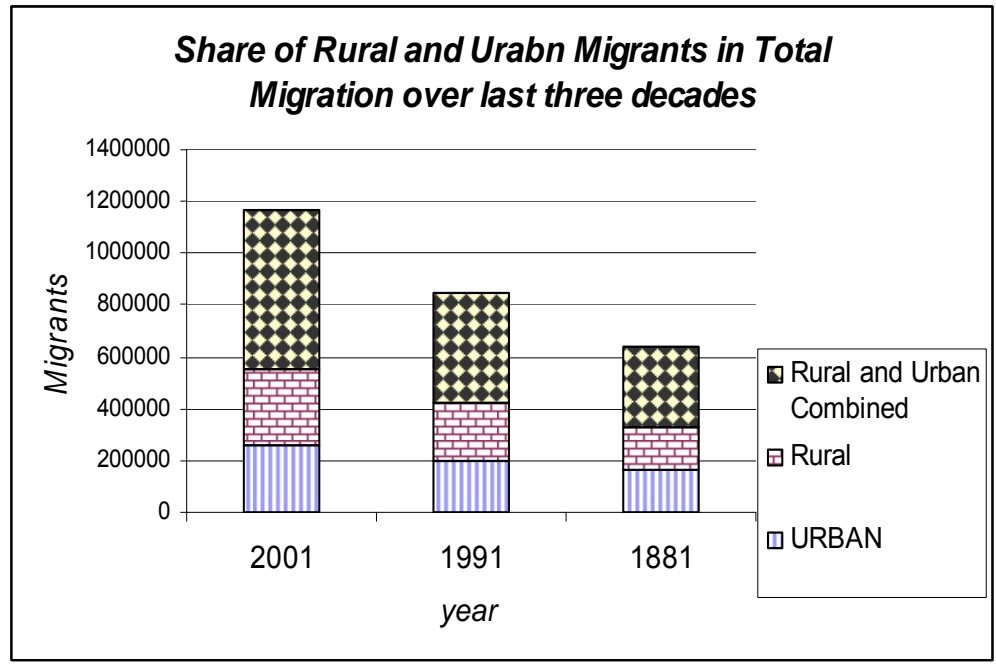

Share of Rural and Urabn Migrants in Total Migration over last three decades

The contribution of rural & urban migrants within a category of reason for migration over last two consecutive decades is tested by calculating the z-values for various category of reason for migration for Rural & Urban areas and significance was tested at 5% level of significance. To test the equality of share of Rural/Urban migrants from any reason of migration over a decade period following hypothesis was set up.

$H_0$ : Share of Rural (or Urban) migrants due to any reason of migration in a decade is equal. ($p_1=p_2$)

Against

$H_1$ : $p_1 \neq p_2$

This is tested for two decadal period 1981-9991 & 1991-2001.

$$Z = \frac{p_1 - p_2}{\sqrt{PQ((1/n_1) + (1/n_2))}}$$

where $P = \dfrac{n_1 p_1 + n_2 p_2}{n_1 + n_2}$ and $Q = 1 - P$

$p_1$ is the share of rural/urban migrants due to any reason at a point of time in total migration,

$p_2$ is the share of rural/urban migrants due to that reason after a decade in total migration

To test this hypothesis, Z-value for equality of proportions of migrants from any reason over a decade is calculated and compared with tabulated value at 5% level of significance for the period 1981-9991 & 1991-2001 for rural and urban migrants separately. Four groups according to share of migrants from any reason of migration over a decade period are formed to analyze the Rural-Urban dynamics of the migrant process.

Group1: Share of migrants from any reason of migration from Rural/Urban area over a decade period (in 1981-991 & 1991-2001) is not equal. Means share of peoples migrating from rural & urban areas for a particular reason of migration differ significantly over the period 1981-991 & 1991-2001. Areas falling under this group shows a change in similar direction (i. e. share of urban & rural migrants for that reason of migration has changed considerably over a decade period) for Rural & Urban migrants in terms of their share in total migration for that reason of migration over a decade period.

Group 2: Share of migrants from any reason of migration from Rural/Urban area over a decade period (in 1981-991 & 1991-2001) is equal. Means share of peoples migrating from rural & urban areas for a particular reason of migration don't differ significantly over the period 1981-991 & 1991-2001. Areas falling under this group don't shows any change (i. e. share of urban & rural migrants for that reason of migration is has not changed over a decade period) for Rural & Urban migrants in terms of their share in total migration for that reason of migration over a decade period.

Group 3: Share of migrants from any reason of migration from Rural area is not equal whereas for migrants from urban areas due to this reason is equal over a decade period (in 1981-991 & 1991-2001). Means share of peoples migrating from Rural areas for a reason of migration differ significantly  whereas share of peoples migrating from Urban areas for this reason of migration don't differ significantly over the period 1981-991 & 1991-2001. Areas falling under this group shows different story as share of Urban migrants for any reason of migration in total migration is not equal though for Rural Migrants it is equal over a decade period.

Group 4: Share of migrants from any reason of migration from Urban area is not equal whereas for migrants from Rural areas due to this reason is equal over a decade period (in 1981-991 & 1991-2001). Means share of peoples migrating from Urban areas for a reason of migration differ significantly whereas share of peoples migrating from Rural areas for this reason of migration don't differ significantly over the period 1981-991 & 1991-2001. Areas falling under this group shows different story as share of Rural migrants for any reason of migration in total migration is not equal though for Urban Migrants it is equal over a decade period.

In Group 1 & 2, migration due to any reason from rural and urban areas is in agreement i.e. share of migrants due to any reason over a decade either is significant or insignificant for both rural and urban migrants. In contrary to this In Group 3 & 4, migration due to any reason from rural and urban areas is not in agreement i.e. share of migrants due to any reason over a decade is significant for urban migrants than it is insignificant for rural migrants or vice-versa.

Z-value for testing hypothesis at 5% level of significance in a group will be as under.

Group 1: $Z_u$ & $Z_r$ >1.96

Group 2: $Z_u$ & $Z_r$ <1.96

Group 3: $Z_u$ >1.96 & $Z_r$ <1.96

Group 4: $Z_u$ <1.96 & $Z_r$ >1.96

Where $Z_u$ and $Z_r$ is the calculated value of Z for migrants due to a reason from Urban & Rural area. The significance of Null hypothesis for all the groups is summarized in table on ensuing page.

| Reason for Migration | Contribution of Rural & Urban Migrants over a decade period is in agreement for any reason of Migration | | | |
|---|---|---|---|---|
| | Duration 1991-2001 | | Duration 1981-1991 | |
| | $Z_u$ & $Z_r$ >1.96 | $Z_u$ & $Z_r$ <1.96 | $Z_u$ & $Z_r$ >1.96 | $Z_u$ & $Z_r$ <1.96 |
| Employment | Total Migration, Elsewhere Jaipur District, Gujrat | Haryana, U.P., Delhi | Total Migration, Elsewhere in Jaipur District, in other Districts, Gujarat, Hrayana, U.P., Punjab, Delhi | |
| Education | Total Migration, Punjab | Gujarat, Haryana, U.P., Delhi | -do- | |
| Marriage | U.P., Punjab, Haryana, Delhi | | Elsewhere in Jaipur District, in other Districts, Gujarat, Hrayana, U.P., Punjab, Delhi | |
| Moved with Family | Total Migration, Elsewhere in Jaipur District, in other Districts, U.P., Punjab, Delhi | | Total Migration, Elsewhere in Jaipur District, in other Districts, Gujarat, Hrayana, U.P., Punjab, Delhi | |
| Reason for Migration | Contribution of Rural & Urban Migrants over a decade period is not agreement for any reason of Migration | | | |
| | Duration 1991-2001 | | Duration 1981-1991 | |
| | $Z_u$ >1.96 & $Z_r$ <1.96 | $Z_u$ <1.96 & $Z_r$ >1.96 | $Z_u$ >1.96 & $Z_r$ <1.96 | $Z_u$ <1.96 & $Z_r$ >1.96 |
| Employment | | in other Districts, Punjab | | |
| Education | | Elsewhere in Jaipur District, in other Districts, | | |
| Marriage | Gujarat | | Total Migration | |
| Moved with Family | Haryana | | | |

It is apparent from this summarization that share of rural & urban migrants in the period 1981 & 1991 differ widely for migrants coming from various places. Especially for migrants coming from other states the share of rural & urban population in yr 1981 & 1991 differ significantly for all the four categories of reason for migration. However this fact was a little bit different in the period of 1991-2001 as migrants coming for education & employment from rural & urban areas of various states don't differ significantly in terms of their share in year 1991 & 2001 in total migration.

Migrants from rural & urban areas due to marriage, employment & education were not in agreement as from some of the areas the proportion of rural migrants in year 1991 & 2001 was significant whereas for urban it was not. Therefore for the duration 1991-2001 migrants from some of the places are not making significant difference in terms of their contribution for some of the reasons to migrate or for rural migrants it is not significant whereas for urban migrants it is significant or vice-versa. This situation was missing in the duration 1981-9991.

**GENDER ISSUES AND DISPARITY:**

Flow of male and female migration governed by different reasons differently and exhibit a different trait over the years. Looking at total in-migration in Jaipur it is found that contribution of males were phenomenal high in the category of people migrating due to employment and education as against the share of female was higher than males in category of persons migrating due to marriages and moving with family. Moreover the fact of male moving faster than females for employment & education and females move faster than male for marriages and moving along family was also evident in short, medium and long distances migration and this gap at the segregated levels was much explicit than the aggregated level. Following hypothesis was formulated to test the gender disparity in migration.

$H_0$ : Share of males (or females) migrants due to any reason for migration in a decade is equal (i.e. $p_1=p_2$)

Against

$H_1$ : $p_1 \neq p_2$

Formula for Z remain same whereas $p_1$ is the share of male/female migrants due to any reason at a point of time in total migration and $p_2$ is the share of male/female migrants due to that reason after a decade in total migration.

To test this hypothesis Z-value for equality of proportions of migrants from any reason over a decade is calculated and compared with tabulated value at 5% level of significance for the period 1981-9991 &

1991-2001 for male and female migrants separately. Four groups according to share of migrants of any reason for migrations over a decade period are formed to analyze the Rural-Urban dynamics of the migrant process.

Z-value for testing hypothesis at 5% level of significance for the four groups will be as under.

    Group 1: $Z_m$ & $Z_f > 1.96$

    Group 2: $Z_m$ & $Z_f < 1.96$

    Group 3: $Z_m > 1.96$ & $Z_f < 1.96$

    Group 4: $Z_m < 1.96$ & $Z_f > 1.96$

Where $Z_m$ and $Z_f$ is the calculated value of Z for male & female migrants due to a reason. The significance of Null hypothesis for all the groups is summarized in table inserted below.

| Reason for Migration | Contribution of Male & Female Migrants over a decade period is in agreement for any reason of Migration | | | |
|---|---|---|---|---|
| | Duration 1991-2001 | | Duration 1981-1991 | |
| | $Z_m$ & $Z_f$ >1.96 | $Z_m$ & $Z_f$ <1.96 | $Z_m$ & $Z_f$ >1.96 | $Z_m$ & $Z_f$ <1.96 |
| Employment | Total Migration, in other Districts | Gujarat, Punjab, Haryana, U.P., Delhi | Total Migration, in other Districts, Gujarat, Punjab, | |
| Education | | Gujarat, Punjab, Haryana, U.P., Delhi | Total Migration | Gujarat, Punjab, Haryana, U.P., Delhi |
| Marriage | Total Migration, | | Total Migration, Elsewhere in Jaipur District, in other Districts, Haryana, Punjab, | |
| Moved with Family | Total Migration, Elsewhere in Jaipur District | | Total Migration, Elsewhere in Jaipur District, in other Districts, Gujarat, U.P, Delhi | |
| Reason for Migration | Contribution of Male & Female Migrants over a decade period is not in agreement for any reason of Migration | | | |
| | Duration 1991-2001 | | Duration 1981-1991 | |
| | $Z_m$ >1.96 & $Z_f$ <1.96 | $Z_m$ <1.96 & $Z_f$ >1.96 | $Z_m$ >1.96 & $Z_f$ <1.96 | $Z_m$ <1.96 & $Z_f$ >1.96 |
| Employment | Elsewhere in Jaipur District in other districts | | Gujarat, Punjab, Haryana, U.P., Delhi | |
| Education | Total Migration, in other Districts | | Elsewhere in Jaipur District, in other Districts | |
| Marriage | | Elsewhere in Jaipur District, in other Districts Gujarat, Punjab, Haryana, U.P., Delhi | | Gujarat |
| Moved with Family | | **in other Districts Gujarat, Punjab, Haryana, U.P., Delhi** | | Punjab, Haryana, U.P., Delhi |

It is evident from the above results that the contribution of male and females in different categories over two decades (1981-91 & 1991-2001) has changed considerably and the disparity is widened. As most of the categories in duration (1981-1991) fall in the group where both $Z_m$ & $Z_f$>1.96 which means proportion of the males & females over a decade was significantly different. In this way male & females for most of the categories were in agreement ($Z_m$ & $Z_f$ >1.96) as both were significant as far as their contribution in total migration over a decade is concerned. Except for the people moving due to education from other states as $Z_m$ & $Z_f$ <1.96 for this category. This means that share male & females migrating due to education from other states in total migration in the year 1981 & 1991 was same and this remained stabilized in year 2001. People migrating due to marriages & moving with family also showed a change in this three decade period as migrating from most of the areas in year 1991 over 1981 exhibited that the share was considerably changed ($Z_m$ & $Z_f$>1.96) whereas in year 2001 over 1991 it showed that it has not changed for males though for females it has changed. Thus people moving under these categories have shown a shift in term of increasing share toward females.

Migrants from different areas exhibit a considerable shift in terms of contribution of males or females in total migration over a period of ten years. However overall migrants say that three categories (employment, marriage & moved with family) followed the same suite as the share of male & female was significant for testing hypothesis for equality of the same over the duration 1981-991 & 1991-2001.

**Summary:**

Contribution of people migrating for education in total migration is on a steep declining as its contribution in total migration has decreased by one third over a two decade period. People migrating due to marriage is showing a phenomenal incremental growth & it is supposed to grow with a faster pace due to decline sex ratio in the city. Migration due to education is having less contribution in total migration and it is going thinner over the years because of education facilities in smaller town and easy accessibility to them in small town. Therefore no longer education is as significant for tempting to migrate as it used to be two decades back. In the coming years this cause of migration will further tend to lose its impact in overall mobility of peoples. People migrating with family is also on a downward trend as people moving with family and due to marriage are together constitute inactive movement as people are not necessarily moving by choice or primarily don't have motive of employment, business or education which itself are related to betterment of life/career.

The share of inactive movements in total migration has came down by 5% over a decade. If this trend continues and the economic progress of the Jaipur indicates that it will attract the people for economic reasons than the share of migrants in working population will grow which in turns contribute for the economy of the City as the share of people moving with family is declining sharply. Migration from urban areas due to marriages is also getting bigger and voluminous in coming decades this will in turn affect the cultural & social structure of the society and a cosmopolitan culture will emerged.

Analysis of trend of the male & females' migration it can be interpreted that share of employment & education from other states to Jaipur is leading to stabilization & it was not found significant for testing the hypothesis of equality of their share over decades. Whereas for people moving with family the share of males is getting stabilized though for females it was growing. If this scenario continue than growing migration of females in this category will, to some extent, be beneficial to the decreasing sex ratio the city.

Short distance migration which consist the in-migration from various parts of the Jaipur district to Jaipur Urban area is one a sharp decline path in terms of its contribution in total migration. It clearly indicates the tendency of migrating to Jaipur urban area is lower down as periphery of Jaipur urban area is also being developed as its suburb. Better connectivity is raising the number of daily commuter and in near future entire district may be developed as a part of Jaipur urban area and a new Jaipur is shaping up. In such a scenario overall migration to Jaipur urban area from the various parts of Jaipur district will lose its relevance.

**References:**

Bhagat, R. B. (992). Components of Urban Growth in India with Reference to Haryana: Findings from Recent Censuses, Centre for Training and Research in Municipal Administration, University of California, Nagarlok, Vol.24, No. 3, 1992, pp. 10-14

Bhattacharya, P. (1988). The informal sector and rural to-urban migration: Some Indian evidence, *Economic and Political Weekly* vol.33(21), 1255-1262

Jain, M.K., Ghosh, M. and Won Bae, K. *(1993)*. Emerging Trends of Urbanization In India: An Analysis Of 1991 Census Results, *Office of Registrar General, New Delhi*

Registrar General of India: Census of India (1971, 1981, 1991, 2001) Migration tables Part V A and B, Rajasthan State.